\documentclass[a4paper,pdflatex,french,12pt]{article}
\usepackage[T1]{fontenc}
\usepackage[latin1]{inputenc}
\usepackage[francais]{babel}
\usepackage{geometry}
\usepackage{amsfonts,amsmath,amssymb,amsthm}
\usepackage{times}
\usepackage{pifont}
\usepackage{graphicx}
\usepackage{subfig}
\usepackage{vmargin,calc}
\usepackage{algorithm,algorithmic}
\usepackage{url}
\usepackage{wrapfig}
\usepackage[all]{xy}

\makeatletter
\setboolean{ALC@noend}{true}
\renewcommand{\ALG@name}{Algorithme}
\makeatother

\def\R{\mathbb R}
\def\Z{\mathbb Z}

\def\upp{\raisebox{-0.3mm}{\includegraphics{dessin_xavier.1}}}
\def\down{\raisebox{-0.3mm}{\includegraphics{dessin_xavier.2}}}
\def\upoudown{\raisebox{-0.3mm}{\includegraphics{dessin_xavier.8}}}

\title{Combinatoire du point de croix}
\author{Sandrine Caruso\footnote{Élève de l'ÉNS Cachan, antenne de 
Bretagne}\ \ et Xavier Caruso\footnote{Chercheur au CNRS, affecté à 
l'Université Indépendante de Moscou}}

\newtheorem{theo}{Théorème}[section]
\newtheorem{lemme}[theo]{Lemme}
\theoremstyle{definition}
\newtheorem{defi}[theo]{Définition}
\theoremstyle{remark}
\newtheorem*{rem}{Remarque}

\begin{document}

\maketitle

\begin{abstract}
Dans cet article, nous expliquons une technique de broderie classique :
le point de croix. Nous nous intéressons ensuite à la question de
minimiser la longueur de fil utilisée pour broder un dessin donné et
résolvons le problème lorsque le dessin est $4$-connexe (notion définie
dans le texte). Nous décrivons également un algorithme qui brode le
dessin avec la quantité minimale de fil attendue.
Enfin, dans une dernière partie, nous étudions plusieurs exemples de
dessins qui ne sont pas $4$-connexes.

\medskip

\noindent
{\it Mots-clés :} géométrie discrète, optimisation, combinatoire
\end{abstract}

\setcounter{tocdepth}{2}
\tableofcontents

\section{Position du problème}

\subsection{Le point de croix}

\begin{figure}[ht]
\centering
\includegraphics{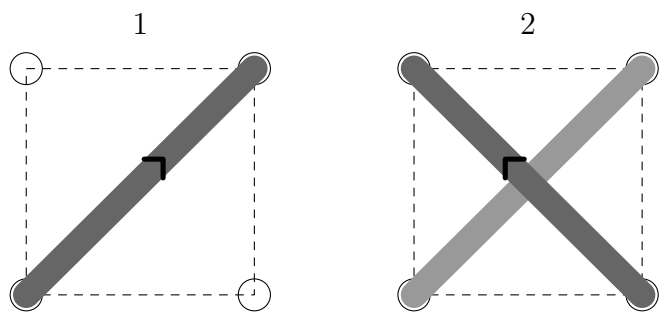}
\caption{Le point de croix}
\label{explic_point}
\end{figure}

Le point de croix est une technique de broderie consistant à reproduire
le dessin d'une grille sur un tissu à trame régulière, à l'aide (comme
son nom d'indique) de points en forme de croix. Pour les lecteurs qui ne
connaîtraient pas le principe de la broderie, voici quelques précisions.
On dispose d'un fil enfilé dans une aiguille, que l'on fait passer à
travers les trous (disposés régulièrement le long d'un quadrillage) du
tissu, alternativement de l'envers vers l'endroit du tissu et de
l'endroit vers l'envers.
Pour broder un point de croix, le fil doit ainsi passer deux fois sur
l'endroit du tissu : une fois selon l'une des diagonales, une fois selon
l'autre diagonale (voir figure~\ref{explic_point} ; la partie de fil
plus foncée représente la dernière diagonale brodée).

\begin{figure}[ht]
\centering
\includegraphics{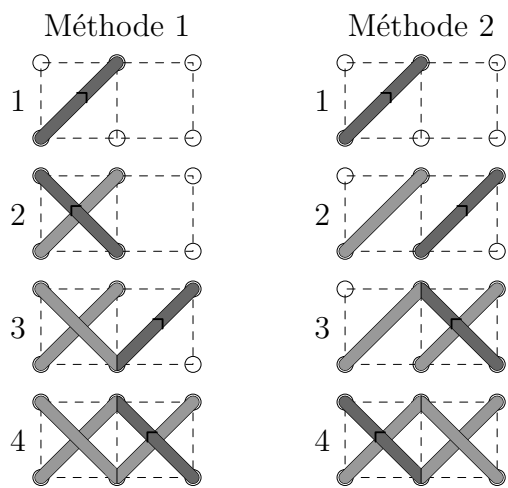}
\caption{Deux méthodes pour broder deux cases}
\label{deux_points}
\end{figure}

Lorsqu'il y a plusieurs cases à broder, l'ordre et le sens dans lequel 
on brode chaque diagonale ne doit obéir qu'à une seule condition : la 
diagonale qui passe par dessus l'autre doit être la 
même\footnote{Sous-entendu : dans la même direction} pour toutes les 
cases. Par exemple, pour broder deux cases côte-à-côte, les deux 
méthodes présentées sur la figure~\ref{deux_points} sont acceptables.

\begin{figure}[ht]
\centering
\includegraphics{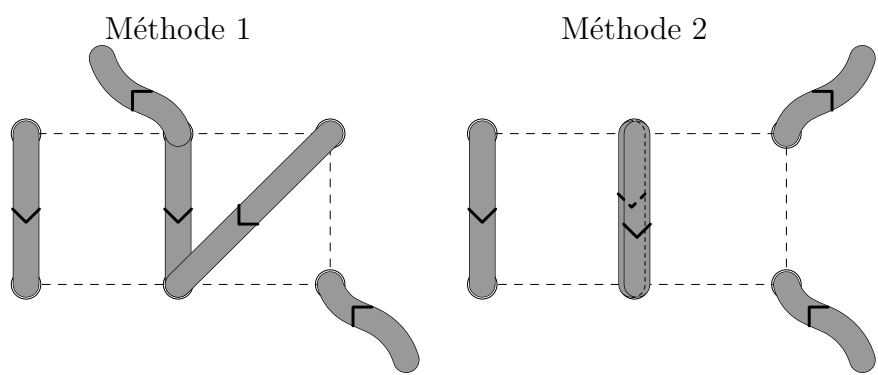}
\caption{Envers du tissu}
\label{deux_points_envers}
\end{figure}

Comme vous pouvez le constater sur la figure, ces deux méthodes donnent
le même motif, et la condition sur les diagonales est respectée.
Cependant, si les deux méthodes produisent le même résultat sur
l'endroit du tissu, en revanche, sur l'envers de celui-ci, les fils ne
sont pas disposés de la même façon, comme on peut le voir sur la
figure~\ref{deux_points_envers}.

Dorénavant, nous allons fixer le choix de la diagonale passant par
dessus l'autre, par exemple, comme dans la figure~\ref{deux_points},
nous supposerons qu'il s'agira de la diagonale joignant le coin en haut
à gauche au coin en bas à droite (si l'on regarde l'endroit du tissu),
que nous appellerons désormais \og diagonale supérieure \fg. Notons que
cela signifie que pour une case donnée, il faudra broder la diagonale
joignant le coin en haut à droite au coin en haut à gauche (que nous
appellerons \og diagonale inférieure \fg) avant l'autre diagonale.

\subsubsection*{Longueur de fil utilisée}

D'un point de vue pratique, on peut avoir envie d'utiliser la méthode
qui demande le moins de fil possible. Faisons quelques remarques à ce
sujet. Nous négligeons la longueur de fil dépassant avant et après la
partie brodée. En outre, si l'on brode le même motif avec plusieurs
méthodes différentes, la longueur de fil présente sur l'endroit du tissu
sera bien entendu la même. En fait, cette longueur ne dépend que du
nombre de cases brodées. Seule est susceptible de varier la longueur de
fil présente sur l'envers du tissu. Par exemple, si l'on choisit comme
unité la longueur du côté d'une maille, pour la méthode 1, la longueur
de fil sur l'envers du tissu est $2 + \sqrt{2}$ tandis que pour la
méthode 2, cette longueur est $3$ (voir la
figure~\ref{deux_points_envers}). Ainsi, la méthode 2 est-elle plus
économe en fil que la méthode 1.

Notons tout d'abord qu'entre deux diagonales que l'on brode
successivement, la longueur minimale sur l'envers est de $1$ (obtenue
dans le cas où le fil passe simplement le long d'un côté de maille). En
effet, on ne peut pas passer l'aiguille deux fois de suite par le même
trou (de l'endroit vers l'envers puis de l'envers vers l'endroit), car
le fil ne serait pas maintenu et ressortirait, lâche, du côté endroit.
Si l'on veut broder $n$ cases, chaque case étant composée de $2$
diagonales, il faut au minimum $2n-1$ unités de fil sur l'envers.

\subsection{Configurations brodables et fortement brodables}
\label{subsec:brodable}

\begin{figure}[h]
\centering
\includegraphics{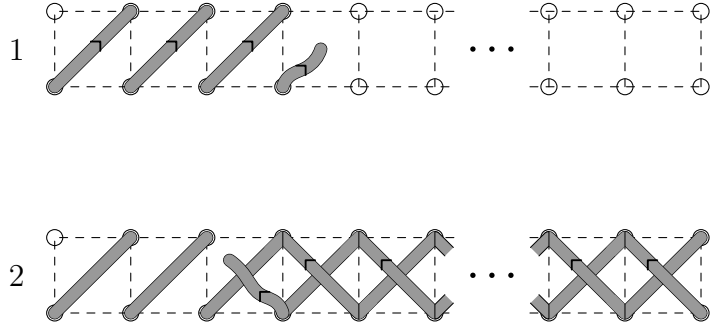}
\caption{Ligne horizontale de cases}
\label{ligne}
\end{figure}

La question que nous nous posons à présent est la suivante : quelles
sont les configurations qu'il est possible de broder avec seulement
$2n-1$ unités de fil sur l'envers ? De telles configurations existent
bel et bien. Un petit manuel de point de croix, par exemple, nous
apprend que c'est le cas notamment si les cases forment une ligne
horizontale (ou verticale) continue ; il suffit de broder d'abord toutes
les diagonales inférieures, puis toutes les diagonales supérieures,
comme sur la figure~\ref{ligne}.

\begin{defi}
On appelle \emph{configuration brodable} une configuration à $n$ cases
qui peuvent être brodées par un fil de longueur $2n (1+\sqrt 2) -1$
(ce qui correspond à une longueur $2n-1$ sur l'envers, c'est-à-dire
la longueur minimale).
\end{defi}

Nous venons de voir que les configurations \og ligne \fg\ et \og colonne
\fg\ sont brodables. En fait, la méthode que nous venons d'expliquer
possède une propriété supplémentaire intéressante : elle est telle que
le point de départ et le point d'arrivée du fil sont deux sommets
adjacents. Autrement dit, si l'on souhaite revenir à notre point de
départ, on a seulement besoin d'une longueur $1$ de fil supplémentaire
(ce qui est bien sûr le minimum). Comme cette notion sera essentielle
dans la suite de notre article, nous dégageons une nouvelle définition.

\begin{defi}
On dit qu'une configuration à $n$ cases est \emph{fortement brodable}
si elle peut être brodée par un fil de longueur $2n(1 + \sqrt 2)$ qui
revient à son point de départ.
\end{defi}

Bien entendu, une configuration fortement brodable est brodable : il
suffit de ne pas faire le dernier point sur l'envers. Nous verrons dans
la partie \ref{sec:escalier} (théorème \ref{theo:escsimple}) que la
réciproque est fausse.

\section{Configurations 4-connexes}

L'objectif de cette partie est de donner une condition suffisante pour
qu'une configuration soit fortement brodable (et donc en particulier
brodable). Nous décrirons ensuite un algorithme qui calcule une façon
de broder une telle configuration.

\subsection{Définitions et énoncé du théorème}

%

Afin de mathématiser notre problème, on modélise le tissu à l'aide du
plan $\R^2$, les trous de la trame étant identifiés aux points du réseau
$\Z^2$. La définition suivante donne un sens précis à quelques termes
déjà employés dans la section précédente. Nous nous tiendrons désormais
à ce vocabulaire.

\begin{defi}[Case, sommet, configuration]
On appelle \emph{sommet} un élément $(n,m) \in \Z^2$. On appelle
\emph{case} un carré de $\R^2$ de la forme $[n,n+1]\times [m,m+1]$ avec
$(n,m) \in \Z^2$. On appelle \emph{configuration} un ensemble fini de
cases.

Si $p = [n,n+1]\times[m,m+1]$ est une case, les quatre sommets $(n,m)$,
$(n,m+1)$, $(n+1,m)$ et $(n+1,m+1)$ sont appelés \emph{sommets de $p$}.
On appelle \emph{sommet d'une configuration} un sommet d'une case de
cette configuration.

Deux sommets sont dit \emph{adjacents} s'ils sont à distance $1$ l'un de
l'autre. On appelle \emph{arête} un segment de droite joignant deux
sommets adjacents. Une \emph{arête d'une case} est une arête joignant
deux sommets de cette case.
\end{defi}


\begin{defi}[$4$-chemin, $8$-chemin]
Soient $p$ et $p'$ deux cases. On appelle \emph{$4$-chemin} reliant $p$
à $p'$ une suite finie de cases $p = p_0, p_1, \ldots, p_n = p'$ telle
que pour tout $i \in \{1, \ldots, n\}$, $p_{i-1}$ et $p_i$ aient au
moins une arête commune.

On appelle \emph{$8$-chemin} reliant $p$ à $p'$ une suite finie de
cases $p = p_0, p_1, \ldots, p_n = p'$ telle que pour tout $i \in \{1,
\ldots, n\}$, $p_{i-1}$ et $p_i$ aient au moins un sommet commun.
\end{defi}


\begin{defi}[$4$-connexité, $8$-connexité]
Une configuration $C$ est dite \emph{$4$-connexe} (resp.
\emph{$8$-connexe}), si pour toutes cases $p$ et $p'$ de $C$, il existe
un $4$-chemin (resp. un $8$-chemin) formé de cases de $C$ reliant $p$ à
$p'$.
\end{defi}
Voici un exemple de figure $4$-connexe et un exemple de figure $8$-connexe 
mais non $4$-connexe.
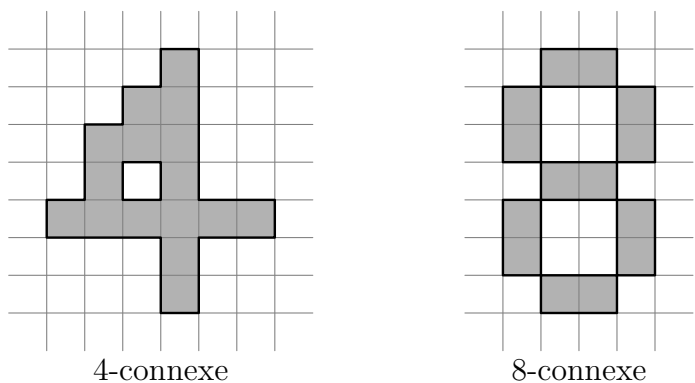
\begin{center}
\includegraphics{dessin.10}
\end{center}

L'objectif de cette section est de démontrer le théorème suivant.

\begin{theo}
Toute configuration $4$-connexe est fortement brodable (et donc
brodable).
\end{theo}

Pour les besoins de la démonstration, définissons encore la notion de
ligne qu'en fait, nous avons déjà rencontré dans l'article.

\begin{defi}
Une configuration $C$ est appelée une \emph{ligne} si elle est 
$4$-connexe et s'il existe un ensemble de la forme $\R \times [m,m+1]$
qui contient toutes ses cases.
\end{defi}

\subsection{Schémas de piquage}

Nous avons vu en \ref{subsec:brodable} qu'une ligne est fortement
brodable. En observant attentivement la méthode qui a été décrite, on se
rend compte que l'aiguille a traversé tous les trous de la ligne
inférieure en allant de l'envers vers l'endroit et tous les trous de la
ligne supérieure en allant de l'endroit vers l'envers, ce que nous
représenterons dorénavant par le schéma de piquage suivant :

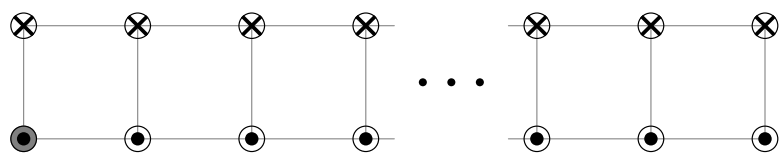
\begin{center}
\includegraphics{dessin_xavier.3}
\end{center}

\noindent

le sommet grisé représentant le point de départ (qui est aussi, 
rappelons-le, le point d'arrivée). Remarquez que la plupart des trous 
ont été traversé deux fois par l'aiguille. Dans notre cas, les deux
passages ont toujours eu lieu dans le même sens, mais c'est à vrai 
dire un peu un hasard et rien n'interdira par la suite d'avoir des 
schémas de piquage où les deux symbôles \upp\ et \down\ apparaissent 
ensemble sur un même trou. On les notera alors l'un à côté de l'autre 
dans un sens arbitraire. Voici une définition précise.

\begin{defi}
Un \emph{schéma de piquage} $\mathcal S$ est la donnée d'une fonction
qui à chaque sommet $s$ à coordonnées entières associe un sous-ensemble
$\mathcal S(s)$ de \{\upp, \down\} et d'un \og point de départ \fg,
c'est-à-dire d'un couple $(s, \text \upoudown)$ où $s$ est un sommet et
$\text \upoudown \in \mathcal S(s)$.

On dit qu'une configuration de $n$ cases est \emph{fortement brodable
selon un schéma $\mathcal S$} s'il est possible de la broder entièrement 
par un fil de longueur $2n(1+\sqrt 2)$ en revenant au point de départ et 
en respectant les indications du schéma.
\end{defi}

\begin{rem}
Si $\mathcal S(s)$ est vide, cela signifie que le sommet $s$ n'est jamais 
piqué. En outre, si une configuration $C$ est fortement brodable selon le 
schéma $\mathcal S$, alors $\mathcal S(s)$ est vide si et seulement si $s$ n'est pas 
un sommet de $C$. En effet, il est évident que tous les sommets de $C$ sont 
piqués. Réciproquement, si $s$ est un sommet piqué, notons 
$p_1,p_2,p_3,p_4$ les quatre cases dont il est sommet.
\begin{itemize}
\item Si l'on pique $s$ de l'envers vers l'endroit, cela signifie que l'on 
s'apprête à broder une diagonale d'un des $p_i$ ($1 \leqslant i \leqslant 4$) 
ou bien qu'on l'a brodée en premier dans le cas où $s$ est le dernier sommet 
piqué (en effet, dans ce cas, $s$ est également le premier sommet piqué, dans 
le même sens, puisque $C$ est fortement brodable).
\item Si l'on pique $s$ de l'endroit vers l'envers, alors on vient de broder une 
diagonale d'un des $p_i$, ou bien on va la broder à la fin dans le cas où $s$ 
est le premier sommet piqué (car, dans ce cas, il est aussi le dernier).
\end{itemize}
Dans tous les cas, un des $p_i$ appartient à la configuration $C$, et donc $s$ est 
un sommet de $C$.
\end{rem}

\subsubsection{Plusieurs schémas de piquage pour broder une ligne}
\label{subsec:ligne}

Nous nous intéressons encore à une configuration de $n$ cases
consécutives. Nous savons déjà qu'une telle configuration est fortement
brodable mais, pour la suite, nous allons besoin d'être un peu plus
précis sur les schémas de piquage associés. C'est l'objet du lemme
suivant.

\begin{lemme}
\label{lem:ligne}
Soit $n$ un entier supérieur ou égal à $2$ et soit $k \in \{2, \ldots,
n\}$. Alors, une ligne de $n$ cases est fortement brodable selon le
schéma suivant :

\begin{center}
\includegraphics{dessin_xavier.4}

\medskip

{\rm Modèle 1 de schéma de piquage}
\end{center}
\end{lemme}

\begin{proof}
On brode les diagonales inférieures des cases $k$ jusqu'à $n$, puis les
diagonales supérieures des cases $n$ jusqu'à $k+1$.

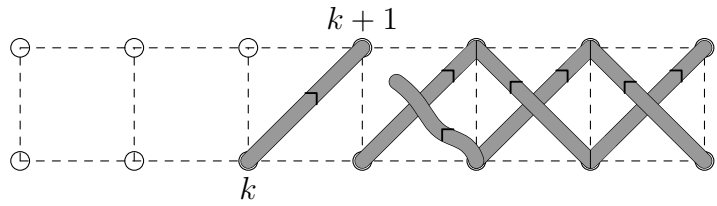
\begin{center}
\includegraphics{dessin.5}
\end{center}

\noindent
L'aiguille vient alors de passer à travers le sommet $k+1$ en haut. On
la passe par le $k$-ième sommet en haut

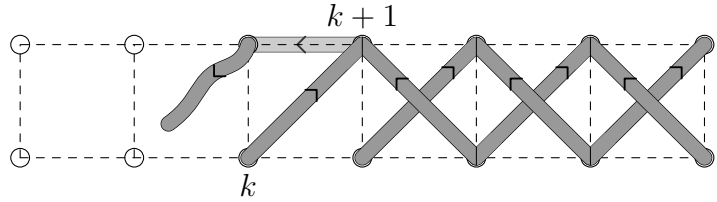
\begin{center}
\includegraphics{dessin.6}
\end{center}

\noindent
de sorte à pouvoir broder les diagonales inférieures des sommets $k-1$
jusqu'à $1$, puis leurs diagonales supérieures. Il ne reste plus qu'à
broder la diagonale supérieure de la case $k$, ce que l'on fait en 
commençant par piquer au point d'abscisse $k+1$ sur la ligne du bas.

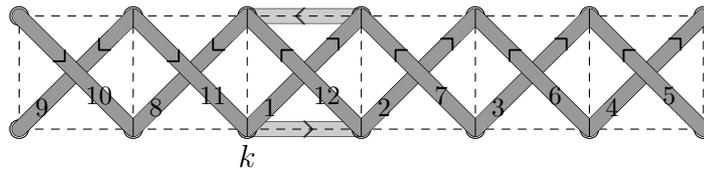
\hfill \includegraphics{dessin.72} \hfill
\end{proof}

\noindent
On aura également besoin de considérer des cas où l'on commence à broder
par en dessus (c'est-à-dire que l'on commence par un point à l'envers).
En déplaçant au début le dernier piquage de l'aiguille dans la méthode
présentée dans la démonstration du lemme, on démontre que la ligne
est aussi fortement brodable selon le schéma que voici :

\begin{center}
\includegraphics{dessin_xavier.5}

\medskip

Modèle 2 de schéma de piquage
\end{center}

\noindent
Enfin, en appliquant une symétrie centrale --- ce qui ne modifie pas le
sens des diagonales --- à la construction que l'on vient de présenter,
on s'aperçoit que la ligne est également fortement brodable selon les
deux schémas suivants :

\begin{center}
\includegraphics{dessin_xavier.6}

\medskip

Modèle 3 de schéma de piquage


\bigskip

\includegraphics{dessin_xavier.7}

\medskip

Modèle 4 de schéma de piquage
\end{center}

\bigskip

Jusqu'à présent, nous avons interdit à $k$ de prendre les valeurs
extrêmes $1$ et $n+1$, c'est-à-dire que nous n'avons pas encore
considéré le cas où le point de départ est à l'extrémité de la ligne à
broder. Dans la situation du lemme \ref{lem:ligne}, on ne peut en fait
par avoir $k = n+1$. En effet, commençant dans le coin en bas à droite,
on serait contraint à broder d'abord la diagonale supérieure du dernier
carré, ce qui pourtant d'après les règles ne peut se faire qu'après
avoir brodé sa diagonale inférieure. Le cas $k = 1$, par contre, ne
conduit pas à une impossibilité mais à un schéma de piquage dégénéré,
déjà bien connu puisque c'est le premier que nous avons rencontré. 
Pour faciliter la lecture, nous le reproduisons ci-dessous.

\begin{center}
\includegraphics{dessin_xavier.9}

\medskip 

{\rm Modèle 1' de schéma de piquage}
\end{center}

\noindent
On remarque que les doubles décorations ont disparu. Comme précédemment,
en déplaçant à la fin le premier point, on obtient le schéma que voici

\begin{center}
\includegraphics{dessin_xavier.10}

\medskip 

{\rm Modèle 2' de schéma de piquage}
\end{center}

\noindent
tandis qu'à l'aide de symétries centrales, on trouve les deux autres
schémas suivants :

\begin{center}
\begin{tabular}{clc}
\includegraphics{dessin_xavier.11} & \hspace{1cm} &
\includegraphics{dessin_xavier.12} \medskip \\
Modèle 3' de schéma de piquage & & Modèle 4' de schéma de piquage
\end{tabular}
\end{center}

\subsection{Démonstration du théorème}

\subsubsection{Recollement des schémas de piquage}

L'intérêt de décorer les configurations par des schémas de piquage
réside dans la propriété de recollement très facile suivante.

\begin{lemme}
\label{lem:recollement}
Soient $C_1$ et $C_2$ deux configurations disjointes fortement brodables
selon les schémas de piquage $\mathcal S_1$ et $\mathcal S_2$
respectivement. Soit $(s, \text \upoudown)$ le point de départ de $\mathcal S_2$.

Si $\text \upoudown \in \mathcal S_1(s)$, alors la configuration $C_1 \cup 
C_2$ est fortement brodable selon le schéma de piquage $P \mapsto 
\mathcal S_1(s) \cup \mathcal S_2(s)$ avec même élément de départ que 
celui de $\mathcal S_1$.
\end{lemme}

\begin{proof}
On fixe une méthode de brodage de $C_1$ (resp. $C_2$) qui respecte le
schéma de piquage $\mathcal S_1$ (resp. $\mathcal S_2$). On commence à
broder $C_1$. D'après l'hypothèse, se faisant, on passera par le sommet
$s$ dans le sens \upoudown. À ce moment, on brode $C_2$. Une fois cela
fait, on est revenu au point de départ en piquant à nouveau dans le sens
\upoudown. On peut ainsi finir de broder $C_1$.

Bien entendu, toutes les conditions sont respectées (tous les points sur
l'envers ont pour longueur $1$ et les diagonales inférieures sont
toujours brodées avant les diagonales supérieures) puisqu'elles
l'étaient déjà pour le brodage de $C_1$ et $C_2$.
\end{proof}

\subsubsection{La démonstration proprement dite}

On fixe à partir de maintenant une configuration $C$ qui est
$4$-connexe. Comme illustré sur la figure ci-après, découpons $C$ en
l'union disjointe de lignes $C_i$ ($1\leq i \leq N$) de sorte que les
cases immédiatement à gauche et à droite de $C_i$ ne soient pas dans
$C$.

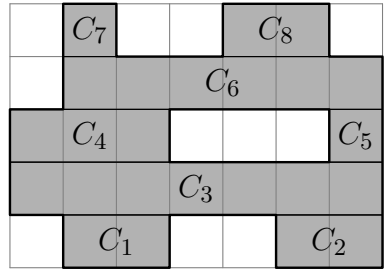
\begin{center}
\includegraphics{dessin.9}
\end{center}

\noindent
Nous allons montrer que pour tout $i \in I = \{1, \ldots, N\}$, on peut
choisir un schéma de piquage $\mathcal S_i$ sur $C_i$ donné par l'un
des modèles vus en \ref{subsec:ligne}, tous ces schémas se
recollant bien. Plus précisément nous allons construire par récurrence :
\begin{itemize}
\item une suite croissante $(I_k)_{1 \leq k \leq N}$ de sous-ensembles
de $\{1, \ldots, N\}$ avec $I_1 = \{1\}$ et $I_k$ de cardinal $k$ pour
tout $k$ (et donc $I_N = I$) ;
\item des schémas de piquage $\mathcal S_i$ sur $C_i$ donnés par l'un
des modèles de \ref{subsec:ligne} tels que le point de départ de
$\mathcal S_1$ soit de la forme $(s, \text \upp)$ et, pour tout 
$k$, la configuration $C' = \cup_{i \in I_k} C_i$ soit fortement
brodable selon le schéma $\mathcal S' : s \mapsto \cup_{i \in I_k}
\mathcal S_i(s)$ ayant pour point de départ celui de $\mathcal S_1$.
\end{itemize}

\medskip

La construction de $I_1$ et de $\mathcal S_1$ ne pose aucun problème :
on prend $I_1 = \{1\}$ comme cela est imposé et on choisit par exemple
le modèle 1'.
Supposons maintenant que $I_k$ soit construit et que les $\mathcal S_i$
pour $i \in I_k$ ont déjà été choisis. Notons $C' = \cup_{i \in I_k}
C_i$ et $\mathcal S'$ le schéma de piquage $s \mapsto \cup_{i \in I_k}
\mathcal S_i(s)$ ayant pour point de départ celui de $\mathcal S_1$.
D'après l'hypothèse, $C'$ est fortement brodable selon $\mathcal S'$.

Soit $p$ une case de $C$ qui n'est pas dans $C'$. Par hypothèse de
$4$-connexité, il existe un $4$-chemin reliant $p$ à une case de $C'$.
Soit $q = p_j$ la dernière case de ce $4$-chemin qui n'est pas dans
$C'$. Appelons $C_{i_0}$ la ligne à laquelle appartient $q$. Comme $q
\not\in C'$, on a nécessairement $i_0 \not\in I_k$. On définit $I_{k+1}
= I_k \cup \{i_0\}$ ; c'est bien un ensemble de cardinal $k+1$.

\begin{wrapfigure}{r}{1.7cm}
\vspace{-0.5cm}

\includegraphics{dessin_xavier.13}
\end{wrapfigure}

On a par ailleurs $q' = p_{j+1} \in C'$, et le fait que $C'$ soit une
union de $C_i$ montre que $q'$ est situé soit au-dessus, soit au-dessous
de $q$. Supposons par exemple que $q$ soit au-dessus de $q'$, l'autre
cas se traitant de manière analogue. Nommons $s$ et $t$ les deux sommets
communs à $q$ et $q'$ comme sur la figure ci-contre. Si $\text{\upp} \in
\mathcal S'(s)$, il suffit de choisir pour $\mathcal S_i$ le modèle 1 ou
1' avec le sommet $s$ pour point de départ. De même, si $\text \down \in
\mathcal S'(t)$, on peut choisir le modèle 4 ou 4' avec $t$ pour point
de départ. Dans les deux cas, le lemme \ref{lem:recollement} montre que
$\cup_{i \in I_{k+1}} C_i = C' \cup C_{i_0}$ est fortement brodable
selon le schéma de piquage voulu. Étant donné que ni $\mathcal S(s)$, ni
$\mathcal S(t)$ ne peut être vide, le seul cas restant est $\mathcal
S'(s) = \{\text \down\}$ et $\mathcal S'(t) = \{\text \upp\}$. Mais, on voit tout de
suite qu'aucun des modèles de schémas de piquage utilisés ne contient à
la suite sur la ligne du haut les deux symbôles \down\ et \upp\ dans cet
ordre. Ainsi ce dernier cas ne peut se produire et la démonstration de
l'hérédité est terminée.

\subsection{Un algorithme pour broder}

La démonstration précédente a l'avantage de fournir avec peu d'effort un
algorithme qui calcule une façon convenable de broder la configuration
qui utilise la quantité minimale de fil. Nous présentons ci-après deux
versions de cet algorithme.

\subsubsection{Version récursive}

La version récursive est légèrement plus facile à comprendre, et c'est
la raison pour laquelle nous commençons par celle-ci.
Dans ce qui suit la lettre {\tt C} désigne une \emph{variable globale}.
Une autre quantité qui doit être considérée comme variable globale est
la position de l'aiguille : quand dans l'algorithme, on dit de piquer en
tel sommet, on déplace \og physiquement \fg\ l'aiguille jusqu'à ce
sommet et elle reste à cette position jusqu'au prochain piquage. 

L'algorithme \ref{algo:appel} est une petite routine qui initialise les
variables et se termine en appelant la fonction {\tt
Broder\_ligne\_haut} (décrite dans l'algorithme \ref{algo:recur}) qui,
couplée à la fonction {\tt Broder\_ligne\_bas}, (décrite au même
endroit) constitue le c\oe ur de l'algorithme.

\begin{algorithm*}[ht!]
\caption{Procédure d'appel}
\label{algo:appel}
\begin{algorithmic}[1]
\STATE {\tt C} $\leftarrow$ la configuration à broder
\STATE {\tt p} $\leftarrow$ une case de {\tt C}
\STATE {\tt s} $\leftarrow$ le sommet en bas à gauche de {\tt p}
\STATE piquer l'aiguille en {\tt s} de l'envers vers l'endroit
\STATE appeler la fonction {\tt Broder\_ligne\_haut}
\end{algorithmic}
\end{algorithm*}

\begin{algorithm*}[ht!]
\caption{Fonction {\tt Broder\_ligne\_haut} (resp. {\tt
Broder\_ligne\_bas})}
\label{algo:recur}
\begin{algorithmic}[1]
\STATE {\tt L} $\leftarrow$ la plus grande ligne horizontale incluse dans
{\tt C} contenant le sommet où se trouve l'aiguille et au-dessus (resp.
au-dessous) de ce sommet
\STATE {\tt liste} $\leftarrow$ la liste ordonnée des sommets à piquer
pour broder {\tt L} à partir de {\tt s} à l'aide d'un des modèles de
\ref{subsec:ligne}
\IF{le brodage n'est pas possible} 
\STATE sortir de la fonction \ENDIF
\STATE retirer de {\tt C} toutes les cases de {\tt L}
\FORALL{sommet {\tt s} dans {\tt liste} parcourue dans l'ordre}
\STATE piquer l'aiguille en {\tt s}
\IF{le sommet {\tt s} est en haut de {\tt L}}
\STATE appeler la fonction {\tt Broder\_ligne\_haut}
\ELSE[le sommet {\tt s} est en bas de {\tt L}]
\STATE appeler la fonction {\tt Broder\_ligne\_bas}
\ENDIF
\ENDFOR
\end{algorithmic}
\end{algorithm*}

\bigskip

Étant donné qu'il n'y a des appels récursifs que lorsque le nombre de
cases dans la variable {\tt C} diminue strictement (et que ceci ne peut
pas se produire une infinité de fois), il est clair que l'algorithme
s'arrête. La correction de l'algorithme, quant à elle, découle de la
preuve du théorème que nous avons donné précédemment. Nous laissons au
lecteur les détails de cette transcription.

On notera finalement que dans l'étape d'initialisation, on choisit une
case quelconque de la configuration $C$. Ceci signifie que, non
seulement, toute configuration $4$-connexe est fortement brodable mais
qu'en outre, il est possible de la broder en commençant par n'importe
laquelle de ses cases.

\subsubsection{Version itérative}

La version itérative, présentée dans l'algorithme \ref{algo:itera}, est
en fait très proche de la version récursive. Elle fonctionne à l'aide
d'une pile dont les éléments sont des couples $(s,f)$ où $s$ est un
sommet et $f$ un élément de l'ensemble $\{ \text{au-dessus},
\text{au-dessous} \}$.

\begin{algorithm*}[ht!]
\caption{Version itérative}
\label{algo:itera}
\begin{algorithmic}[1]
\STATE {\tt C} $\leftarrow$ la configuration à broder
\STATE {\tt p} $\leftarrow$ une case de {\tt C}
\STATE {\tt s} $\leftarrow$ le sommet en bas à gauche de {\tt p}
\STATE créer une pile vide {\tt pile}
\STATE empiler ({\tt s}, au-dessus) sur {\tt pile}
\WHILE{{\tt pile} est non vide}
\STATE dépiler {\tt pile} 
\STATE ({\tt s}, {\tt f}) $\leftarrow$ l'élément dépilé
\STATE piquer l'aiguille en {\tt s}
\STATE {\tt L} $\leftarrow$ la plus grande ligne horizontale
incluse dans {\tt C} contenant le sommet où se trouve l'aiguille et située
{\tt f} de ce sommet
\STATE {\tt liste} $\leftarrow$ la liste ordonnée des sommets à piquer 
pour broder {\tt L} à partir de {\tt s} à l'aide d'un des modèles de
\ref{subsec:ligne} (on convient que {\tt liste} est vide si le brodage
n'est pas possible)
\IF{{\tt liste} est non vide}
\STATE retirer de {\tt C} toutes les cases de {\tt L}
\ENDIF
\FORALL{sommet {\tt s} dans {\tt liste} parcourue dans le sens inverse}
\IF{le sommet {\tt s} est en haut de {\tt L}}
\STATE empiler ({\tt s}, au-dessus) sur {\tt pile}
\ELSE[le sommet {\tt s} est en bas de {\tt L}]
\STATE empiler ({\tt s}, au-dessous) sur {\tt pile}
\ENDIF
\ENDFOR
\ENDWHILE
\end{algorithmic}
\end{algorithm*}

Nous laissons à nouveau l'exercice au lecteur de montrer que
l'algorithme termine et a bien le comportement voulu. Cela est plutôt
facile lorsque l'on a bien compris le fonctionnement de la pile, et
notamment fait le lien entre la pile qui apparaît dans la version
itérative et la pile des appels récursifs de la version récursive.

\subsection{L'algorithme et les configurations non $4$-connexes}

Nous avons pour l'instant examiné le comportement de l'algorithme
lorsqu'on l'appelle avec une configuration $4$-connexe, mais il fait
encore sens de l'appeller avec une configuration ne vérifiant pas cette
propriété. Que se passe-t-il dans ce sens? Remarquons déjà qu'il est
clair que l'algorithme s'arrête encore (l'argument donné précédemment
n'utilisait pas la $4$-connexité) et qu'il brode au moins la composante
$4$-connexe contenant la ligne $C_1$. Il se peut cependant qu'il en
brode plus comme le montre l'exemple très simple de la configuration
suivante (représentée en gris)

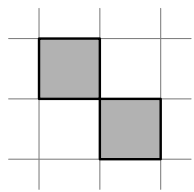
\begin{center}
\includegraphics{dessin_xavier.15}
\end{center}

\noindent
pour laquelle les deux cases sont brodées si l'on part du coin en bas à
gauche de la case du bas. En examinant d'un peu plus près la situation,
on se rend compte que les seuls recollements non $4$-connexes entre deux
des modèles présentés en \ref{subsec:ligne} sont les quatre suivants :

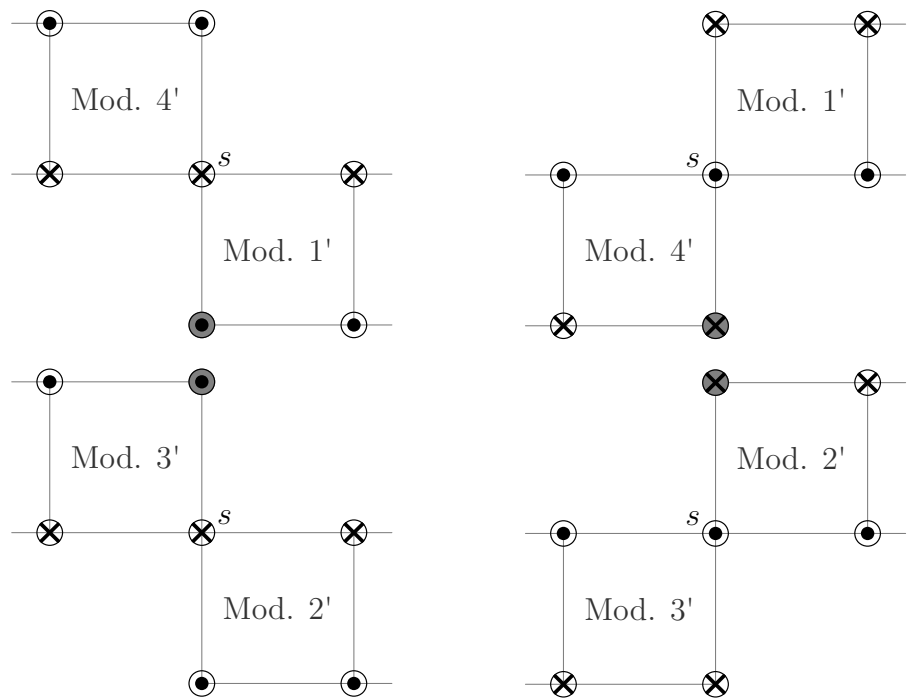
\begin{center}
\raisebox{-0.5\height}{\includegraphics{dessin_xavier.16}}
\hspace{1.5cm}
\raisebox{-0.5\height}{\includegraphics{dessin_xavier.18}}

\bigskip

\raisebox{-0.5\height}{\includegraphics{dessin_xavier.17}}
\hspace{1.5cm}
\raisebox{-0.5\height}{\includegraphics{dessin_xavier.19}}
\end{center}

\noindent 
où $s$ est le sommet de recollement. Ainsi si, dans l'algorithme, on
supprime l'appel récursif après 
\begin{itemize} 
\item le piquage dans le coin en haut à gauche lors de l'exécution du
modèle 1',
\item le piquage dans le coin en bas à gauche lors de l'exécution du
modèle 2',
\item le piquage dans le coin en bas à droite lors de l'exécution du
modèle 3', et
\item le piquage dans le coin en haut à droite lors de l'exécution du
modèle 4'
\end{itemize} on obtient un
programme qui brode \emph{exactement} la composante $4$-connexe de $C$
contenant la ligne $C_1$.

\medskip

Si vous souhaitez voir l'algorithme à l'\oe uvre, rendez-vous sur
la page 

\begin{center}
\url{http://boumbo.toonywood.org/sandrine/pageperso/pcroix/} 
\end{center}

\noindent
où vous pourrez voir se broder devant vos yeux des configurations
préenregistrés ainsi que toutes celles que vos dessinerez puis
proposerez.
La version de l'algorithme utilisée est la version itérative qui prend
en compte la modification que nous venons de discuter pour ne broder
qu'une composante $4$-connexe. En réalité, le programme ne s'arrête pas
après la première composante $4$-connexe, mais continue jusqu'à avoir
bordé toute la configuration proposée en changeant de fil --- et de
couleur --- après chaque composante $4$-connexe.

\section{Exemples et contre-exemples}

Jusqu'à présent, nous n'avons étudié que les configurations $4$-connexes
mais les questions de brodabilité et forte brodabilité ne se posent pas
uniquement dans ce cadre restreint. On a notamment envie maintenant
d'étudier le cas des configurations $8$-connexes. La situation semble
alors bien plus complexe, et c'est ce que aimerions illustrer dans
cette dernière partie à l'aide d'exemples et de contre-exemples.

\subsection{Les configurations \og escalier \fg}
\label{sec:escalier}

Nous étudions, dans cette partie, deux familles d'exemples (que nous
appelerons des \emph{escaliers}) de configurations $8$-connexes qui ne
sont pas $4$-connexes. 

\subsubsection{Un premier lemme bien utile}

Intéressons-nous pour commencer à un type particulier de configurations,
à savoir celles qui sont en accord avec la description suivante :

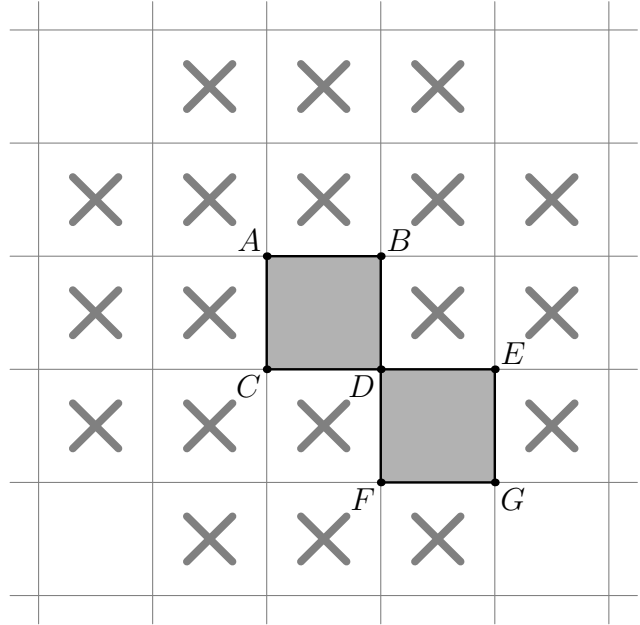
\begin{center}
\includegraphics{dessin_xavier.21}
\end{center}

\noindent
On lit cette description en convenant qu'une case grisée impose la
présence de la case en question dans la configuration, alors qu'une
case barrée impose son absence. Les cases laissées blanches, quant à
elles, n'imposent aucune contrainte.

On considère à partir de maintenant une configuration $C$ de la forme
précédente. On note $n$ le nombre de cases de $C$. On suppose en outre
que $C$ est brodable, et on fixe une manière de la broder qui respecte
les règles que nous nous sommes fixées. Ceci nous permet de numéroter
les diagonales des cases de la configuration : celle qui porte le numéro
$1$ est la première que la méthode de brodage choisie nous dit de
broder, celle qui porte le numéro $2$ est la deuxième et ainsi de suite
jusqu'au numéro $2n$.

\begin{lemme}
\label{lem:escalier}
On se place dans la situation qui vient d'être décrite, et on note $i$
le numéro de la diagonale $AD$. Alors $i \in \{2, 2n\}$. 

De plus, si $i = 2$, alors la diagonale $BC$ porte le numéro $1$, et
la diagonale $EF$ porte le numéro $3$. Si, au contraire, $i = 2n$, 
alors la diagonale $BC$ porte le numéro $2n-1$, et la diagonale $DG$
porte le numéro $2n-2$.
\end{lemme}

\begin{proof}
Supposons par l'absurde que $i \not\in \{2, 2n\}$. Comme il est clair 
que l'on ne peut pas non plus avoir $i = 1$ (puisqu'il faut broder $BC$
avant $AD$), on peut parler des diagonales numérotées $i-2$, $i-1$ et
$i+1$. Il est facile de se convaincre que les diagonales $i-1$ et $i+1$,
étant voisines de $i$, sont à choisir parmi $BC$ et $EF$. Mais, par
ailleurs, les contraintes nous imposent de broder $BC$ avant $AD$. La
seule solution restante est donc que la diagonale $i-1$ soit $BC$ et que
la diagonale $i+1$ soit $EF$. Mais alors la diagonale $i-2$ (qui est
voisine de $i-1$) est nécessairement $DG$. Il en résulte que $DG$ est
brodée avant $EF$, ce qui n'est pas possible. On a ainsi obtenu une
contradiction, et la première partie du lemme est démontrée.

Supposons $i = 2$. Comme la diagonale $BC$ doit être brodée avant $AD$,
elle est nécessairement brodée en premier. D'autre part, après $BC$ on
ne peut broder que $AD$ ou $EF$. Comme $AD$ est déjà brodé, c'est $EF$ 
qui porte le numéro $3$.

Supposons maintenant $i = 2n$. Les deux seuls candidats pour porter le
numéro $2n-1$ sont alors $BC$ et $EF$, mais $EF$ est a écarter car il
doit être brodé après $DG$. Ainsi $BC$ porte bien le numéro $2n-1$ et
il suit rapidement que $DG$ est numéroté $2n-2$.
\end{proof}

Par symétrie centrale, on obtient le résultat suivant : si $C$ est
une configuration brodable qui est de la forme

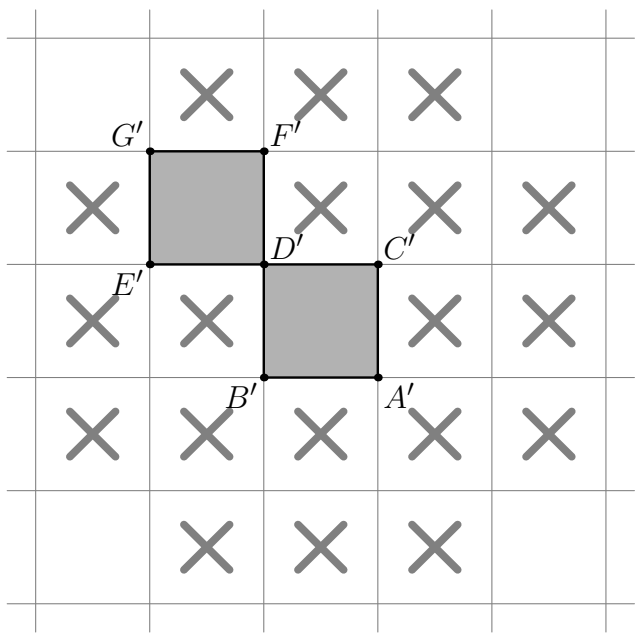
\begin{center}
\includegraphics{dessin_xavier.22}
\end{center}

\noindent
alors le numéro $j$ de la diagonale $A'D'$ appartient à l'ensemble
$\{2, 2n\}$, \emph{etc.}

\subsubsection{L'escalier simple}

Soit $n$ un entier supérieur ou égal à $2$. Considérons la configuration
suivante $E_n$ à $n$ cases.

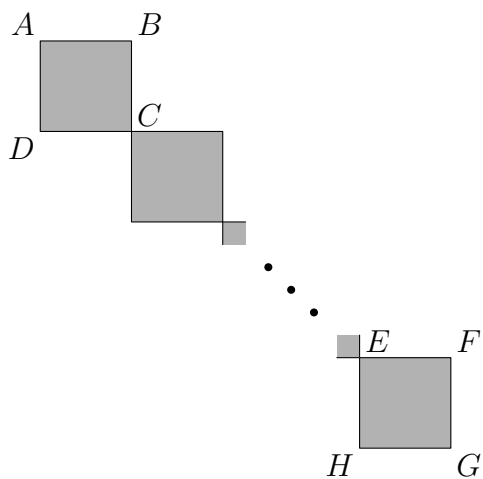
\begin{center}
\includegraphics{dessin_xavier.23}
\end{center}

\noindent
On comprend aisément d'où vient le nom d'\emph{escalier}. On travaillera
dans toute la suite de l'article avec des escaliers qui descendent, mais
les mêmes résultats et les mêmes preuves sont valables pour des
escaliers qui montent.

\begin{theo}
\label{theo:escsimple}
(1) Pour tout $n$, la configuration $E_n$ est brodable.

(2) La configuration $E_n$ est fortement brodable si, et seulement si 
$n = 2$.
\end{theo}

\begin{proof}
L'assertion (1) est facile. On brode les cases successivement en
descendant l'escalier en piquant l'aiguille pour chacune d'elle
\begin{itemize}
\item d'abord dans le coin en bas à gauche,
\item ensuite dans le coin en haut à droite,
\item ensuite dans le coin en haut à gauche,
\item et enfin dans le coin en bas à droite.
\end{itemize}

\medskip

On passe maintenant à la démonstration de (2). Si $n = 2$, voici une
solution pour réaliser le brodage avec les conditions voulues (les
numéros sur les sommets indiquent l'ordre dans lequel l'aiguille les
traverse).

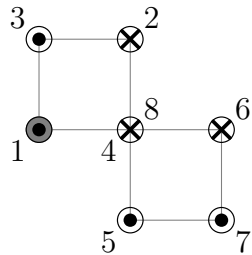
\begin{center}
\includegraphics{dessin_xavier.24}
\end{center}

\noindent
Supposons maintenant $n > 2$ et appelons $ABCD$ et $EFGH$ les cases
respectivement situées en haut et en bas de l'escalier comme cela a déjà
été fait sur la première illustration. Supposons que l'escalier soit
fortement brodable et fixons une manière de le broder qui respecte les
contraintes correspondantes. Si $i$ désigne le numéro de la diagonale
$AC$, le lemme \ref{lem:escalier} nous dit que $i \in \{2, 2n\}$. De
même, le numéro $j$ de la diagonale $EG$ est aussi dans $\{2, 2n\}$.
Quitte à appliquer une symétrie centrale, on peut supposer que $i = 2$
et $j = 2n$. Mais alors, le lemme \ref{lem:escalier} à nouveau nous
apprend que la diagonale $BD$ est bordée en premier. Or ceci est
incompatible avec le fait que $EG$ soit bordée en dernier car ni $B$ ni
$D$ n'est adjacent à $E$ ou $G$ (on rappelle que, par l'hypothèse de
forte brodabilité, le brodage est supposé se terminer au même point
que celui où il a commencé).
\end{proof}

\subsubsection{L'escalier avec palier}

Soient $g$, $p$ et $d$ des entiers naturels non nuls. On considère la
configuration suivante $E_{g,p,d}$

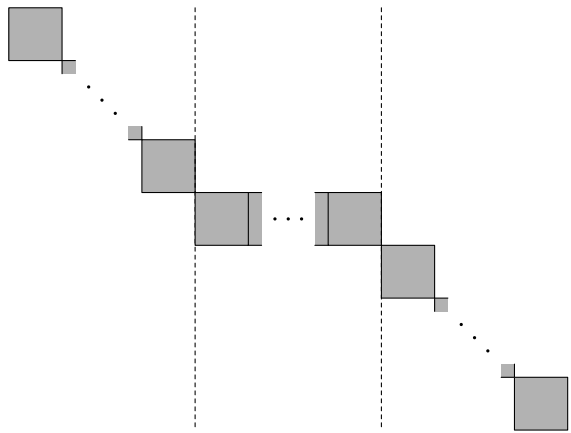
\begin{center}
\includegraphics{dessin_xavier.25}
\end{center}

\noindent
où il y a $g$ cases dans la partie de gauche (le premier escalier), $p$
cases dans la partie centrale (le palier) et $d$ cases dans la partie de
droite (le deuxième escalier). Pour plus de clarté, les trois parties
sont séparées par des traits en pointillés sur le dessin précédent.

Le cas $p = 1$ correspond à un escalier simple où le palier se réduit 
à une marche normale. Ce cas a déjà été étudié précédemment, et c'est pourquoi
nous l'excluons à partir de maintenant.

\begin{theo}
\label{theo:escpalier}
Soient $g$, $p$ et $d$ des entiers naturels non nuls avec $p \geq 2$.
Alors

\begin{itemize}
\item si $g = d = 1$, la configuration $E_{g,p,d}$ est fortement brodable
\item si $g = 1$ et $d > 1$, la configuration $E_{g,p,d}$ est brodable,
mais non fortement brodable
\item si $d = 1$ et $g > 1$, la configuration $E_{g,p,d}$ est brodable,
mais non fortement brodable
\item si $g > 1$ et $d > 1$, la configuration $E_{g,p,d}$ n'est jamais
fortement brodable, et elle est brodable si, et seulement si $p$ est
impair.
\end{itemize}
\end{theo}

Le théorème précédent montre que la combinatoire de l'\og escalier avec
palier \fg\ n'est pas vraiment simple et laisse présager qu'il risque
d'être difficile de trouver une condition nécessaire et suffisante
facilement exploitable pour caractériser les configurations brodables
(resp. fortement brodables) parmi les configurations $8$-connexes. En
tout cas, les auteurs n'ont, à ce jour, pas de réponse satisfaisante
à cette question.

\subsubsection{Démonstration du théorème}

\paragraph{Les cas limites}

On commence par supposer $d = 1$. On souhaite tout d'abord montrer que
$E_{g,p,1}$ est toujours brodable, et pour cela il suffit d'exhiber une
façon convenable de broder cette configuration. Pour éviter de
multiplier encore les notations, on la donne ci-dessous simplement dans
le cas particulier $g = 2$, $p = 3$ en laissant au lecteur l'exercice de
faire les adaptations nécessaires pour le cas général.

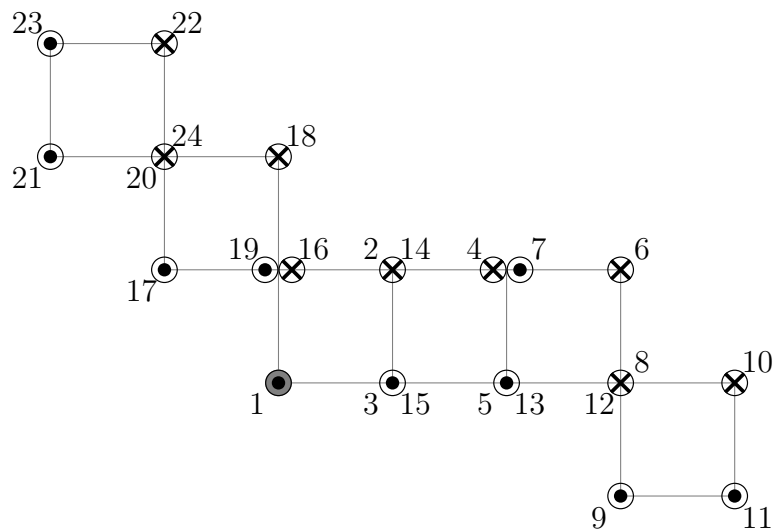
\begin{center}
\includegraphics{dessin_xavier.26}
\end{center}

\noindent
Dans le cas où $d = 1$, une construction similaire montre que
$E_{1,p,1}$ est fortement brodable. Il ne reste donc plus qu'à démontrer
que $E_{g,p,1}$ n'est pas fortement brodable pour $g > 1$. Pour cela,
concentrons-nous sur les deux premières marches en haut à gauche et
nommons leurs sommets comme sur la figure ci-après.

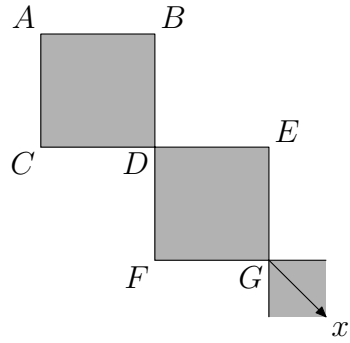
\begin{center}
\includegraphics{dessin_xavier.27}
\end{center}

\noindent
D'après le lemme \ref{lem:escalier}, s'il existe une manière convenable
de broder $E_{g,p,1}$, le numéro d'ordre qu'elle attribue à la diagonale
$AD$ est soit $2$, soit $2n$ (où $n = g + p + 1$ est le nombre total de
cases). Si c'est $2$, d'après le même lemme, les diagonales $BC$ et $EF$
ont respectivement pour numéro $1$ et $3$. On en déduit que la diagonale
$DG$ porte le numéro $4$. En effet, partant de $E$ ou $F$, on ne peut
aller sur l'envers qu'en $D$ ou $G$, et il n'est pas possible de broder
tout de suite la diagonale $Gx$ car il s'agit d'une diagonale supérieure
et que la diagonale inférieure correspondante n'a pas encore été brodée.
Mais alors, il est clair que l'on ne pourra jamais revenir à notre point
de départ puisque celui-ci (qui est $B$ ou $C$) n'est adjacent qu'à $A$
et $D$ et que toutes les diagonales arrivant à ces sommets ont déjà été
brodées.
Un raisonnement analogue conduit également à une contradiction dans le
cas où la diagonale $AD$ est brodée en dernier.

\bigskip

Le cas $g = 1$ se déduit de ce qui vient d'être fait par application
d'une symétrie centrale. Il ne reste donc plus qu'à traiter le cas où
$g$ et $d$ sont tous les deux $> 1$, ce que nous supposons à partir de
maintenant. Remarquons alors que l'argument que nous avons utilisé dans
la démonstration du theorème \ref{theo:escsimple}.(2) s'applique encore
pour montrer que la configuration $E_{g,p,d}$ n'est pas fortement
brodable. Il ne reste donc plus qu'à démontrer qu'elle est brodable si,
et seulement si $p$ est impair.

Supposons tout d'abord que $E_{g,p,d}$ soit fortement brodable et fixons
un brodage qui respecte les conditions correspondantes. Comme cela a
été expliqué précédemment, ce brodage détermine une numérotation des
diagonales des cases de $E_{g,p,d}$.

\paragraph{Coloriage des diagonales}

Commençons par colorier les sommets à l'aide de deux couleurs que l'on
alterne comme sur un échiquier. De façon formelle, la couleur du sommet
de coordonnées $(i,j)$ est donc la classe modulo $2$ de $i+j$. Mais
plutôt que de parler de classe modulo $2$, nous conviendrons dans la
suite que les couleurs que nous utilisons sont simplement le noir et le
blanc.

Une diagonale dans une case relie toujours deux sommets de même couleur,
et on convient alors de donner à cette diagonale cette couleur commune.
Du fait que deux sommets adjacents sont de couleur différente, on déduit
que deux diagonales dans $E_{g,p,d}$ qui sont consécutives (\emph{i.e.}
qui portent des numéros consécutifs) sont de couleur différente.

\paragraph{Ordre de parcours de l'escalier}

Notons $G$ l'ensemble des cases de $E_{g,p,d}$ qui sont dans le premier
escalier à gauche ; par définition, c'est un ensemble de cardinal $g$
dont on note $G_1, \ldots, G_g$ les éléments en convenant que les cases
$G_1, \ldots, G_g$ apparaissent dans l'ordre quand on descend
l'escalier. On définit de même les ensembles $P = \{P_1, \ldots, P_p\}$
et $D = \{D_1, \ldots, D_d\}$.
D'après le lemme \ref{lem:escalier}, les numéros des diagonales
supérieures de $G_1$ et $D_d$ sont $2$ et $2n$ (où $n = g+p+d$ est le
nombre total de cases), et quitte à faire une
symétrie centrale on peut supposer que les numéros sont attribués dans
cet ordre. La deuxième partie du lemme \ref{lem:escalier} assure alors
que la diagonale inférieure de $G_1$ porte le numéro $1$.

\begin{lemme}
Notons $A_i \in \{G, P, D\}$ la partie qui contient la case contenant
la diagonale numéro $i$. Alors
\begin{itemize}
\item pour $1 \leq i \leq 2g$, on a $A_i = G$ ;
\item pour $2g+1 \leq i \leq 2(g+p)$, on a $A_i = P$ ;
\item pour $2(g+p)+1 \leq i \leq 2(g+p+d) = 2n$, on a $A_i = D$.
\end{itemize}
\end{lemme}

\begin{proof}
On sait déjà que $A_1 = G$ et $A_{2n} = D$. Par ailleurs, il est clair
que, pour tout $i$, on ne peut avoir ni $A_i = G$ et $A_{i+1} = D$, ni
$A_i = D$ et $A_{i+1} = G$ (les cases de $G$ et $D$ sont trop
éloignées). Pour conclure, il suffit donc de montrer qu'il existe au
plus deux indices $i$ tels que $\{A_i, A_{i+1}\} = \{G, P\}$ et, de
même, au plus deux indices $j$ tels que $\{A_j, A_{j+1}\} = \{P, D\}$.

À partir de maintenant, on se concentre sur la preuve de ses assertions.
On traite même en fait uniquement la première, la seconde étant
totalement similaire. Remarquons que si $i$ est tel que $A_i = G$ et
$A_{i+1} = P$, alors la diagonale numérotée $i$ est forcément l'une des
diagonales de $G_g$. De même, si $A_i = P$ et $A_{i+1} = G$, alors la
diagonale numérotée $i+1$ est forcément l'une des diagonales de $G_g$.
Ainsi, à chaque $i$ tel que $\{A_i, A_{i+1}\} = \{G, P\}$, on peut
associer une diagonale de $G_g$ notée $d(i)$. Comme $G_g$ possède deux
diagonales, il suffit pour conclure de montrer que cette association est
injective.

Pour cela, raisonnons par l'absurde en considérant deux indices $i$ et
$j$ avec $i < j$ et $d(i) = d(j)$. Il est alors clair, au vu de la
définition, que $d(i) = d(j)$ doit être la diagonale numérotée $i+1 =
j$. Cela signifie que $A_i = P$, $A_{i+1} = G$ et $A_{i+2} = P$. Mais
alors, les diagonales numérotées $i$ et $i+2$ sont des diagonales de
$P_1$, et donc forcément ses deux diagonales. Mais par ailleurs, puisque
leurs numéros diffèrent de $2$ (qui est un nombre pair), elles doivent
être de même couleur, ce qui n'est manifestement pas le cas. On a donc
obtenu une contradiction d'où résulte l'injectivite annoncée puis le
lemme.
\end{proof}

\paragraph{Étude au niveau du palier}

Le lemme précédent nous dit que le brodage de l'escalier avec palier se
découpe en trois parties : on brode d'abord complètement $G$, puis on
passe à la partie $P$ que l'on brode complètement avant de broder
finalement $D$. Examinons de plus près la façon dont on brode le palier
$P$. Nécessairement, puisque l'on vient de $G$, on entre dans $P$ par le
coin à haut à gauche (sommet $s$) juste après avoir brodé la diagonale
supérieure de $G_g$. De même, on sort de $P$ par le coin en bas à droite
(sommet $t$) juste avant de se mettre en position pour broder la
diagonale inférieure de $D_1$. En résumé, on a le schéma de piquage
partiel que voici :

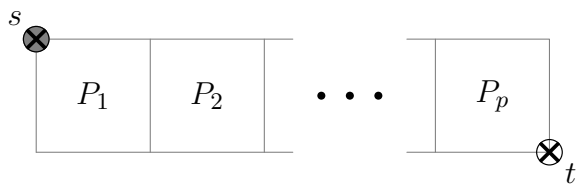
\begin{center}
\includegraphics{dessin_xavier.28}
\end{center}

\noindent
Convenons, pour fixer les idées, que $s$ est colorié en noir. La
première diagonale brodée dans $P$, c'est-à-dire la diagonale numérotée
$2g+1$ est donc blanche. Il en résulte que la dernière diagonale brodée
dans $P$ qui porte le numéro $2g+2d$ est noire puisque la différence
$(2g+2d)-(2g+1) = 2d-1$ est un nombre impair. Or cette dernière
diagonale contient le sommet $t$ ; celui-ci est donc également noir,
c'est-à-dire de la même couleur que $s$. On en déduit que $p$ est impair
comme voulu.

\paragraph{La réciproque}

On suppose désormais que $d$ est un nombre impair, et on souhaite
exhiber une façon de broder $E_{g,d,p}$ répondant aux contraintes
usuelles. Voici une façon de faire. On commence par broder les cases de
$G$ (avec les mêmes notations que précédemment) comme cela est expliqué
dans la preuve du théorème \ref{theo:escsimple}.(1). À l'issue de cela,
le fil sort par le coin en bas à droite de $G_g$ dans la sens \down. On
brode ensuite les cases de $P$ deux par deux selon le schéma suivant :

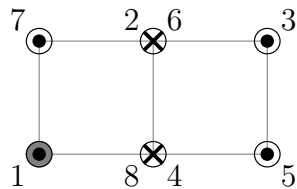
\begin{center}
\includegraphics{dessin_xavier.29}
\end{center}

\noindent
Le fil vient alors de traverser le coin en bas à gauche de $P_{p-1}$
dans le sens \down. Finalement, on utilise à nouveau la méthode de la
preuve du théorème \ref{theo:escsimple}.(1) pour broder ce qui reste
(c'est-à-dire la partie $D$ à laquelle est ajoutée la case $P_p$).

\subsection{Récapitulatif : un diagramme d'implications}

Le diagramme suivant récapitule les implications que nous avons
démontrés entre les principales notions définies dans cet article.

$$
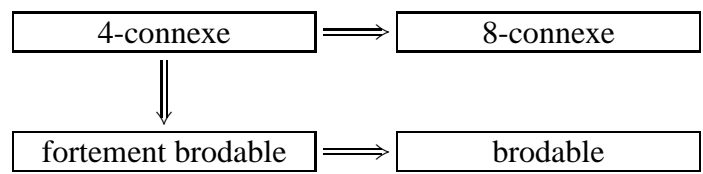
\xymatrix{
\framebox[4cm]{4\text{-connexe}} \ar@{=>}[r] \ar@{=>}[d] & 
\framebox[4cm]{8\text{-connexe}} \\
\framebox[4cm]{\text{fortement brodable}} \ar@{=>}[r] & 
\framebox[4cm]{\text{brodable}}}$$

\medskip

\noindent
Soulignons que les implications qui ne sont pas notées sur le précédent
diagramme sont toutes fausses :
\begin{itemize}
\item le théorème \ref{theo:escsimple} montre qu'il existe des
configurations $8$-connexes qui sont brodables sans être fortement
brodables ;
\item le théorème \ref{theo:escpalier} montre qu'il existe des
configurations $8$-connexes qui ne sont pas brodables ;
\item la configuration formée de deux cases sur une même ligne séparée
par une unique case n'est certainement pas $8$-connexe mais est pourtant
fortement brodable :

\medskip

\begin{center}
\includegraphics{dessin_xavier.30}
\end{center}
\end{itemize}

\end{document}